# Physics-Guided Spectral Parametric Reduced-Order Modeling for Transient Prediction of Controlled Dynamical Systems

Ao Zhang [a], Tian Zhang [a], Antonio Cammi [b], Xiang Wang [a,*]

[a] College of Nuclear Science and Technology, Harbin Engineering University, Harbin, Heilongjiang, 150001, China

[b] Department of Energy, Politecnico di Milano, Milano 20133, Italy

E-mail address: zhangao@hrbeu.edu.cn, tian.zhang@hrbeu.edu.cn, antonio.cammi@polimi.it, xiang.wang@hrbeu.edu.cn

## ABSTRACT

Efficient parametric transient prediction at unseen parameter values and under new operating conditions remains a challenge in scientific computing, uncertainty quantification, design exploration, and control analysis, because repeated high-fidelity simulations are often computationally prohibitive. Existing data-driven surrogates and parametric reduced-order models often perform well within sampled parameter ranges but lose reliability beyond available conditions. This study proposes a physics-guided spectral parametric reduced-order modeling framework for controlled dynamical systems. Parameter-dependent reduced spectral operators are first identified from transient snapshots using Dynamic Mode Decomposition with control, separating intrinsic dynamics from external control effects. After physics-guided parameter transformation, aligned spectral quantities and reduced operator components are propagated across parameter conditions using Secondary Dynamic Mode Decomposition; linear and radial basis function regressions are retained as baselines. Baseline regularization and nondimensional time mapping improve cross-parameter robustness. Validation considers a mechanical transmission system, a Helium-Xenon closed Brayton cycle, and a Kármán vortex street, covering linear transient, nonlinear transient, and nonlinear periodic dynamics. For the mechanical and Brayton systems, system-level multivariable transient responses are predicted at unseen parameter values and under new operating conditions with relative norm errors below 1%. For the vortex-street system, nondimensional time mapping preserves dominant vortex-shedding structures and transient dynamics across Reynolds numbers. Further analyses compare the proposed framework with an LSTM surrogate and assess extrapolation confidence using an auxiliary error-prediction model. Overall, the framework extends parametric reduced-order modeling to extrapolative transient prediction, thereby broadening data-driven reduced-order modeling from in-domain approximation to predictive analysis of complex controlled dynamical systems under unseen conditions.

[*] Corresponding author

**NOMENCLATURE**

| | | | |
|---|---|---|---|
| $A$ | Reduced state transition operator | $t$ | Physical time, s |
| $B$ | Reduced control mapping operator | $U$ | Inlet velocity |
| $b$ | Mode amplitudes | $U_r$ | Orthonormal basis |
| $D$ | Cylinder diameter, m | $X$ | State snapshot matrix |
| $d$ | Damping matrix | $x$ | System state |
| DMD | Dynamic mode decomposition | *Greek symbols* | |
| DMDc | Dynamic mode decomposition with control | $\beta/\beta^{+}$ | Control inputs / Measured auxiliary inputs |
| $G$ | Identified augmented operator ($[A,\ B]$) | $\mu/\mu'$ | Parameters / transformed coordinate |
| $J$ | Moment of inertia, kg·m² | $\Phi$ | Eigenvectors / dynamic modes |
| $K$ | Stiffness matrix | $\Omega$ | Eigenvalues |
| $k$ | Gear ratio | $\Sigma_r$ | Matrix of singular values |
| LBM | Lattice Boltzmann method | $\tau$ | Nondimensional time coordinate |
| $M$ | Inertia matrix | *Subscripts* | |
| pROM | Parametric reduced-order model | $m$ | Motor |
| $Re$ | Reynolds number | $i$ | Intermediate |
| $r$ | Rank | $t$ | Temporal |
| RBF | Radial basis function | $p$ | Parametric |
| $St$ | Strouhal number | $new$ | New parameter values or control inputs |
| SVD | Singular value decomposition | $raw$ | Unregularized / absolute response data |
| $T$ | Temperature, K | $ss$ | Steady offsets |

# 1. Introduction

## 1.1 Background

Parameterized dynamical systems subject to external control arise widely in scientific and engineering applications, where inputs and loads are actively regulated to maintain desired performance, stability, and safety under changing operating conditions. Accurate transient prediction is therefore essential for design exploration, uncertainty quantification, control assessment, and safety analysis. However, the coupled dependence of transient behavior on system parameters, external control inputs, and operating conditions often necessitates repeated high-fidelity simulations across numerous configurations, leading to substantial computational cost. Data-driven surrogate models provide a non-intrusive alternative by accelerating transient prediction without repeatedly solving complex governing equations (Ljung et al., 2020). Nevertheless, most existing methods remain strongly tied to the sampled dataset and are primarily effective for reconstructing or interpolating observed response patterns. When parameter values or operating conditions extend beyond the sampled domain, reliable prediction requires the underlying low-order dynamics to be transferred across conditions rather than merely fitted to available trajectories. Therefore, developing a framework that can transfer coherent reduced dynamics across parameter values and operating scenarios is crucial for efficient parametric transient prediction of controlled dynamical systems.

## 1.2 Requirements for Parametric Transient Prediction

Parametric transient prediction becomes necessary when system responses vary substantially with design parameter values and operating conditions. Existing studies have mainly addressed such sensitivity through static multi-objective optimization under nominal operating conditions (Wang et al., 2018). However, Qin et al. (2025) on $sCO_2$ power cycles and Ferdoues et al. (2017) on wind turbines both indicated that design optimization focused only on nominal design points can lead to unsatisfactory system performance under multi-condition operation. Accordingly, they relied on dynamic system simulation models or transient CFD simulations, respectively, to predict and optimize off-design operation and multi-condition transient performance. Nevertheless, for transient-prediction tasks requiring repeated evaluations, extensive dynamic simulations or high-fidelity numerical simulations can introduce substantial computational cost during parameter studies, control assessment, and scenario evaluation. This motivates the development of efficient and reliable methods for parametric transient prediction under varying parameter values and operating conditions.

### 1.3 State of the Art in Parametric Transient Prediction Methods

Machine learning-based surrogates have gained increasing attention for their ability to learn temporal input-output mappings and rapidly infer transient responses at lower computational cost after training (Brunton et al., 2020). However, their black-box nature and heavy data dependence often limit interpretability and transferability, thereby weakening prediction reliability under parameter extrapolation (Lusch et al., 2018). For instance, relevant to the case studies presented in this work, recurrent neural network (RNN) models predicting Kármán vortex streets under fixed boundary settings (Xie et al., 2024), together with long short-term memory (LSTM) models applied to Brayton-cycle transients (Zhang et al., 2026) and polymer electrolyte membrane fuel-cell (PEMFC) voltage responses under load transients (Elkholy et al., 2025), were primarily validated within training-space coverage or similar transient scenarios, effectively functioning as interpolators rather than extrapolators. Similarly, in the domain of Neural Operators, Zhu et al. (2023) quantitatively showed that DeepONets applied to parametric partial differential equations can suffer large errors under out-of-distribution extrapolation, and recovering validity typically requires post-training fine-tuning with physics-based constraints.

Given these challenges, Parametric Reduced-Order Modeling (pROM; Benner et al., 2015) offers an interpretable, operator-based alternative by regressing low-rank operators over the parameter space instead of learning opaque weights, enabling rapid yet credible predictions across different parameter conditions. Existing comparisons further show that reduced-order identification methods can achieve competitive accuracy with lower computational cost than black-box machine-learning surrogates in certain dynamical-system forecasting tasks, including metabolic oscillations (Wüstner et al., 2025), electric loads (Muñoz et al., 2024), and short-term traffic flows (Yu et al., 2021), supporting their suitability for limited-sample parametric generalization.

### 1.4 Limitations of Existing Parametric Reduced-Order Modeling

However, pROM applications to controlled dynamical systems remain limited, particularly for transient prediction under unseen parameter values and new operating conditions. As the foundation for parametric modeling, data-driven pROM typically leverages Proper Orthogonal Decomposition (POD) or Dynamic Mode Decomposition (DMD; Schmid, 2010), favoring the latter for time-series prediction as it yields a reduced set of spatiotemporal modes and an explicit linear model for their temporal evolution (Brunton and Kutz, 2019); however, standard DMD is restricted to autonomous evolution and lacks a mechanism to separate external inputs from internal dynamics (Schmid, 2022; Zhang et al., 2025a). While Dynamic Mode Decomposition with control (DMDc; Proctor et al., 2016) addresses this limitation, parametric frameworks centered on DMDc remain scarce, leaving a gap for controlled system analysis.

Moreover, existing pROM methodologies for parametric adaptation can be broadly divided into global and local strategies (Fig. 1). Global models are constructed from stacked snapshots across parameter conditions (Sayadi et al., 2015), but they do not explicitly represent the coherent variation of reduced dynamics with parameter. Local pROMs, by contrast, usually operate in an interpolative regime within the sampled parameter domain (Riva et al., 2026), where reduced operators are predicted through numerical fitting (Du et al., 2026; Huhn et al., 2023), including linear and radial basis function (RBF) interpolation (Chen et al., 2024), or through black-box ML surrogates (Min et al., 2026), including BPNN-based surrogates (Mo et al., 2025; Kang et al., 2025) and LSTM-based surrogates (Zhang et al., 2025b; Zhou et al., 2023). In these approaches, parameter variation is often treated as a pure regression problem, while the coherent dynamical structures governing parameter changes remain insufficiently captured, which substantially limits their extrapolation capability in prior work. Therefore, existing pROM studies remain limited in transient prediction for controlled dynamical systems under different parameter conditions, especially in the physically coherent generalization of low-order dynamics across the parameter space.

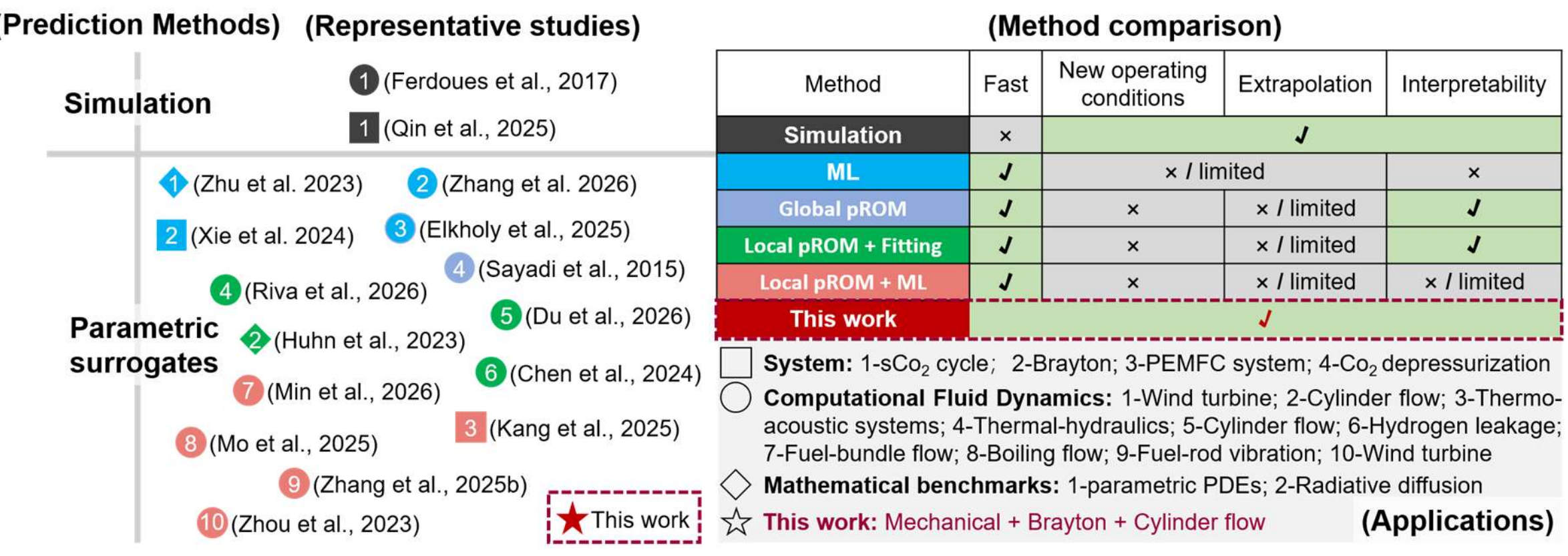

(a) Comparison of existing transient prediction methods and positioning of this work.

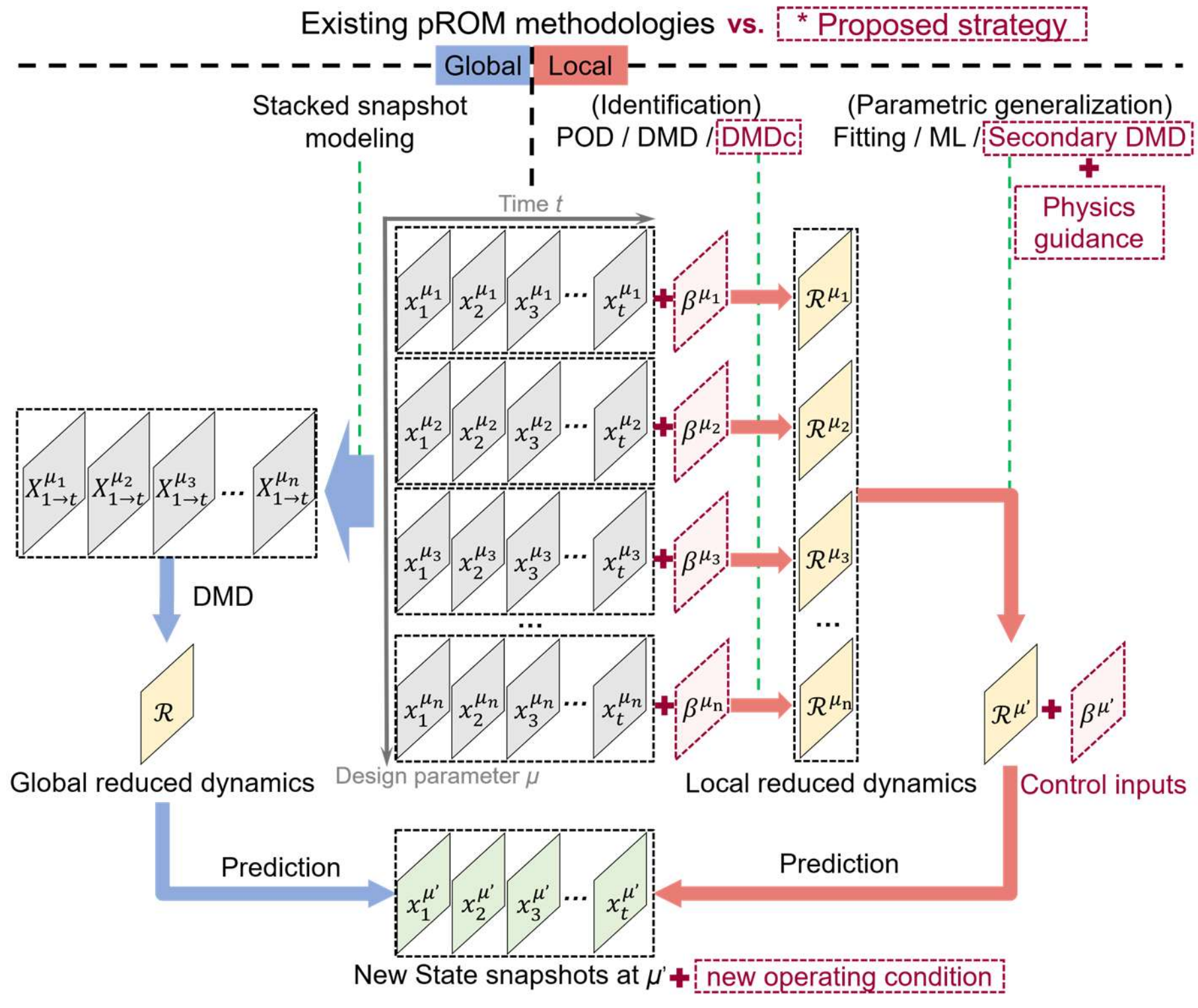

(b) pROMs for parametric generalization.

**Fig. 1.** Review of existing approaches and novelty positioning of this work.

## 1.5 Contributions and Organization

To address the demand for efficient and interpretable transient prediction beyond sampled parameter values and operating conditions, this paper proposes a physics-guided spectral pROM framework. The main methodological advance lies in extending DMDc-based control-separated reduced-order modeling to parametric transient prediction, where Secondary DMD represents parameter-induced variations in aligned spectral quantities and reduced operator components as coherent evolving dynamics. In particular, the work emphasizes extrapolative prediction, for which target reduced operators are no longer bracketed by neighboring sampled conditions and must therefore be reconstructed from spectral evolution beyond the observed parameter interval. The framework is evaluated on three representative systems spanning linear transient, nonlinear transient, and nonlinear periodic dynamics.

We detail this methodology, together with representative case studies with distinct dynamical characteristics (Section 2). Subsequently, we validate the framework for unseen parameter values and operating conditions across these cases, demonstrating robust predictive performance through comparative and sensitivity analyses (Section 3). Finally, a summary of the work and future perspectives is provided (Section 4).

## 2. Theory and Modeling

### 2.1 Fundamentals of Dynamic Mode Decomposition with Control

DMDc identifies reduced-order operators from temporal data, enabling both reconstruction of original transients and prediction under new control inputs. To extract an interpretable low-order model from high-dimensional controlled data, the standard DMDc (Proctor et al., 2016) is formulated here in a reduced-order setting that can be implemented in either discrete or continuous time. Truncated singular value decomposition (SVD) of the state snapshot matrix, $X \approx U_r \Sigma_r V_r^*$, yields an orthonormal basis $U_r$ that denoises the data and maps the high-dimensional state $x$ to reduced coordinates $x_r = U_r^* x$. The reduced dynamics are represented in discrete time and continuous time, respectively, as

$$x_{r,k+1} = A x_{r,k} + B\beta_k \quad , \quad \dot{x}_r(t) = A x_r(t) + B\beta(t) \tag{1}$$

where $A$ and $B$ denote the reduced state transition and control operators.

The reduced state and control data are assembled into the augmented matrix $Y$, and the operator $G = [A, B]$ is identified by least squares as

$$G = \arg\min_{\mathrm{G}} \parallel Z_r - GY \parallel_F^2 = Z_r\, Y^{\dagger} \tag{2}$$

where $Z_r = X_r$ denotes the one-step shifted reduced snapshot matrix in discrete time, and $Z_r = \dot{X}_r$ denotes the reduced derivative snapshot matrix in continuous time. Separating the state transition operator $A$ from the control mapping operator $B$ enables the spectral decomposition, $A\Phi = \Phi\Omega$, where the eigenvectors $\Phi$, eigenvalues $\Omega$, and mode amplitudes $b$ characterize the dominant growth rates, frequencies, and dynamic modes of the system, serving as the primary quantities for subsequent parametric generalization (Riva et al., 2026). For prediction on nonuniform grids, the continuous-time model is advanced under a zero-order hold assumption, with state propagation through the matrix exponential $e^{A\Delta t}$ and the corresponding integrated input term, whereas in the discrete-time setting, prediction follows directly from the identified one-step evolution operator. Importantly, the linear form of DMDc does not imply globally linear controlled dynamics, since its DMD-based formulation can approximate dominant nonlinear temporal behavior through a finite-dimensional data-driven linear representation, which supports its use for nonlinear systems (Schmid, 2010; Proctor et al., 2016).

### 2.2 Proposed Physics-Guided Spectral pROM Framework

The proposed physics-guided spectral pROM framework follows the workflow illustrated in Fig. 2. System spectral dynamics are first categorized (Section 2.2.1) to identify representative dynamical archetypes that guide the modeling strategy. Subsequently, baseline regularization and DMDc-based spectral identification are performed (Section 2.2.2), whereby reduced operators including the state transition operator $A$ and the control mapping operator $B$ are identified along with the basis $U_r$, from which spectral constituents ($\Phi,\ \Omega,\ b$) are extracted to ensure modal consistency and preserve an interpretable representation of the reduced dynamics. Finally, physics-guided regression and generalization (Section 2.2.3) are employed to generalize spectral quantities for new parameters $\mu_{new}$ and control inputs $\beta_{new}$, enabling prediction of the transient response $X_{new}$. Overall, the framework enables efficient prediction at unseen parameters by reusing spectrally aligned reduced operators learned from limited samples, without repeated simulations or case-specific surrogate construction.

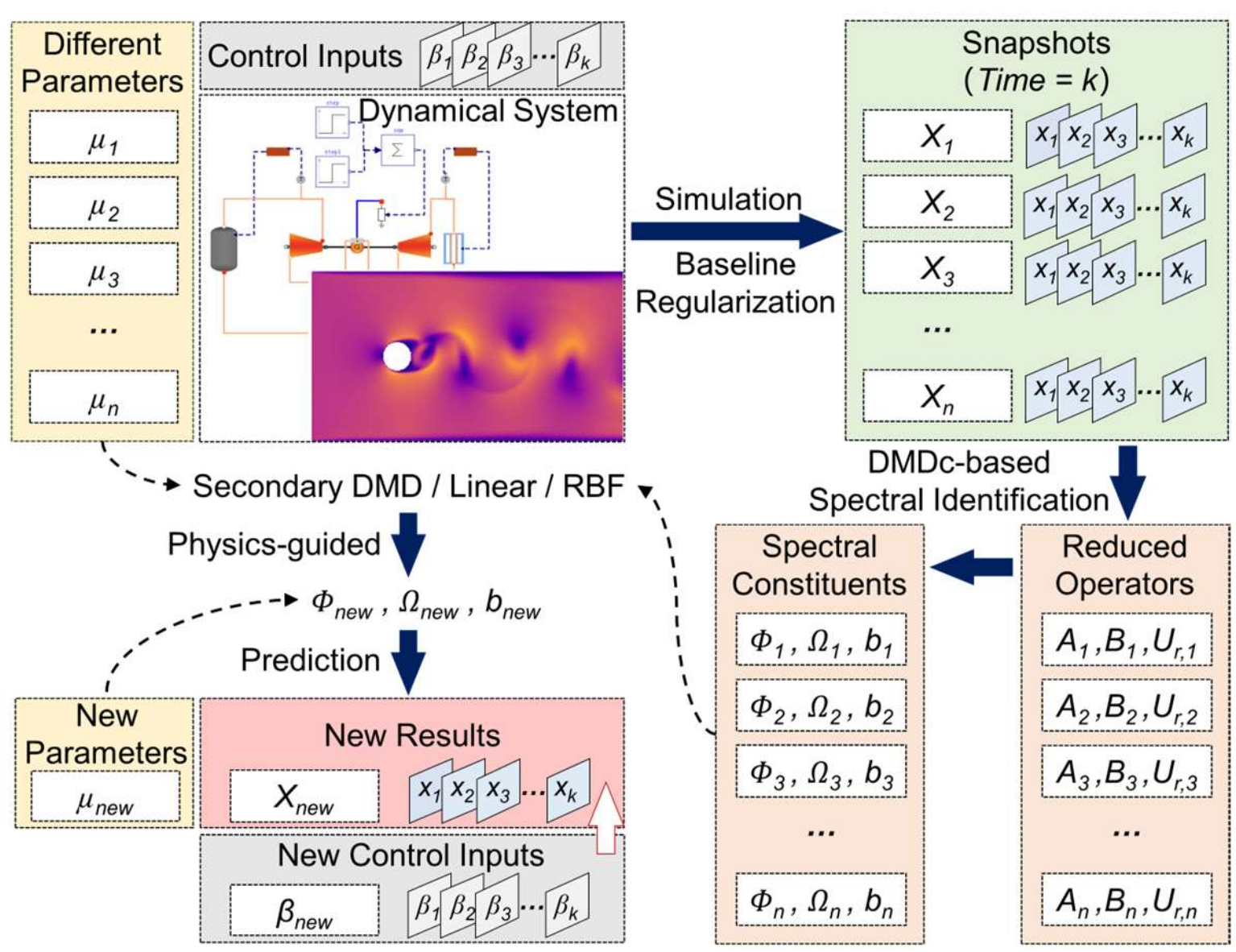

**Fig. 2.** Physics-guided spectral pROM framework.

### 2.2.1 Dynamical Archetypes

To evaluate the proposed framework under representative dynamical behaviors, the systems considered in this study are grouped into three representative dynamical archetypes based on their behaviors: linear transient systems, nonlinear transient systems, and nonlinear periodic systems.

- Linear transient systems obey the superposition principle, with rates of change that are linear functions of the state, thereby serving as baseline cases for validating the basic accuracy of spectral identification when controlled responses remain relatively simple under changing operating conditions.

- Nonlinear transient systems violate the superposition principle, with rates of change governed by nonlinear functions of the state, imposing greater demands on the proposed framework for spectral identification and prediction.
- Nonlinear periodic systems are characterized by trajectories that converge to stable limit cycles and serve as test cases for evaluating the consistency of frequency tracking and modal alignment achieved by the proposed framework under changing operating conditions.

These archetypes cover the main dynamical behaviors considered in this work and provide a clear basis for the subsequent spectral identification and transient prediction.

### 2.2.2 Baseline Regularization and Spectral Identification

Baseline regularization is first introduced so that the snapshot matrix $X$ used for identification represents the transient response relative to the local operating condition, rather than the absolute response containing steady offsets. In the present framework, the regularized state and control snapshots are constructed as

$$X = X_{\mathrm{raw}} - X_{\mathrm{ss}}, \qquad \beta = \beta_{\mathrm{raw}} - \beta_{\mathrm{ss}} \tag{3}$$

where $X_{\mathrm{ss}}$ and $\beta_{\mathrm{ss}}$ denote the steady offsets, which are regressed separately and then used to reconstruct the full response. This treatment reduces the influence of operating-point shifts on subsequent regression and allows the identified reduced operators to focus on the intrinsic transient dynamics.

For nonlinear periodic systems, an additional temporal coordinate transformation is introduced to account for parameter-dependent variations in the characteristic time scale. Since Eq. (2) enforces a single constant reduced operator $G = [A, B]$, significant frequency variation during a transient makes the dynamics in physical time nonuniform, thereby motivating the introduction of a transformed nondimensional time coordinate

$$\tau = \tau(t) \tag{4}$$

under which the periodic transient is recast onto a more uniform dynamical scale. In this way, the full process can be represented more coherently by one reduced operator, while avoiding spectral blurring caused by direct identification in physical time.

Following baseline regularization, spectral identification is performed in reduced coordinates, in discrete time for fixed-interval responses, and in continuous time on the transformed coordinate for nonlinear periodic systems.

The resulting spectral quantities are aligned across training conditions through eigenvalue matching, phase correction, and normalization to preserve modal correspondence. In parallel, the control mapping operator is treated in factorized form,

$$B = U_B \Sigma_B V_B^* \tag{5}$$

allowing its dominant factors to be aligned and generalized separately. When local reduced bases are used, they are aligned across training conditions and reconstructed at the target parameter so that the generalized reduced operators can be applied consistently in the local reduced coordinate system. This stage yields the spectrally consistent reduced quantities $(A, B, U_r, \Phi, \Omega, b)$ for the subsequent regression and generalization stage.

### 2.2.3 Physics-Guided Regression and Generalization

After spectral identification, the aligned reduced quantities are generalized in parameter space through a physics-guided regression strategy, in which physical insight is used to define a transformed coordinate

$$\mu' = f(\mu) \tag{6}$$

that makes their variation along the parameter axis more nearly linear, rather than imposing governing-equation constraints. Accordingly, the ordered training samples can be treated as a pseudo-temporal sequence, allowing Secondary DMD to model the parameter-domain evolution of the ordered reduced quantities through a propagation operator rather than pointwise fitting. Here, Secondary DMD denotes a parameter-domain DMD regression on aligned spectral quantities and factorized operator components, rather than a second temporal decomposition of the original state snapshots. Linear and RBF regressions are used as baseline regressors for comparison.

When a simple parameter mapping is insufficient to regularize all variations, additional nonlinear observables are explicitly incorporated into $\beta$ as measured auxiliary inputs $\beta^+$, so that their nonlinear effects are absorbed by the control term rather than by the state-transition operator itself. Such auxiliary inputs are treated as available measured or prescribed closed-loop signals in the present controlled-response prediction setting. In this way, $A$ can more cleanly represent the intrinsic reduced evolution, while $B$ maps the prescribed and auxiliary inputs into the reduced-state evolution. The generalized spectral and operator components are then recombined to construct the target reduced model, which is advanced using the prediction procedure in Section 2.1 to obtain the transient response at $\mu_{new}$ under $\beta_{new}$.

### 2.3 Validation Cases

#### 2.3.1 Linear Transient System: Mechanical Transmission

Mechanical transmission systems are common components for mechanical energy transfer in controlled dynamical systems. Motor-load inertia matching is a typical industrial design problem for such systems (Wang et al., 2022), because the inertia ratio influences transient torque transmission and speed response. Accordingly, the motor inertia $J_m$ is treated as the design parameter defining different parameter conditions, while the driving torque applied at the motor side serves as the control input. This setting represents a motor-load inertia-matching scenario for parametric transient prediction.

The simulation model (Fig. 3) is implemented in OpenModelica and verified against the standard rotational-mechanics example in the Modelica Standard Library (Fritzson et al., 2019). It exhibits linear torsional dynamics under torque excitation, governed by rotational balance, elastic shaft coupling, viscous damping, and ideal gear kinematics. Accordingly, the governing equations admit the second-order form

$$M\ddot{q} + d\dot{q} + Kq = B_c\beta \tag{7}$$

where $q$ collects the generalized angular coordinates, $M$, $d$, and $K$ denote the inertia, damping, and stiffness matrices, respectively, and $B_c$ is the input distribution matrix. After conversion to first-order state space, the continuous-time operator contains the blocks $M^{-1}K$, $M^{-1}d$, and $M^{-1}B_c$, indicating that inertia enters the dynamics through its inverse. In the present mechanical transmission system, reflection of the motor inertia through the ideal gear yields a dominant motor-dependent inertial contribution of the form $J_i + k^2 J_m$ on the intermediate side, whose inverse is well approximated over the parameter range considered by a term proportional to $1/J_m$, thereby providing the physical basis for the parameter transformation adopted in Eq. (6). For reduced-order identification, the state snapshots comprise the angular displacements and velocities of the three inertias together with selected internal elastic and damping variables, while the driving torque serves as the control input.

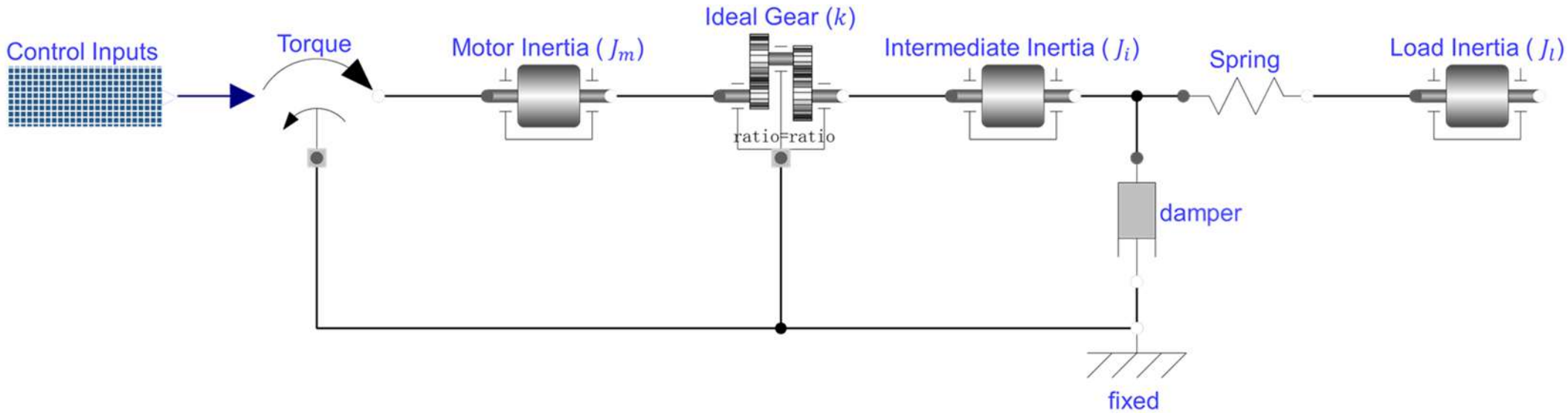

**Fig. 3.** Modelica-based model of the mechanical transmission.

### 2.3.2 Nonlinear Transient System: Closed Brayton Cycle

Energy-conversion systems are representative controlled dynamical systems in the energy domain, and the closed Brayton cycle is a typical example involving strong thermodynamic coupling and active load regulation. In this case, the minimum cycle temperature $T_c$ is treated as the system parameter defining different parameter conditions, while the load-demand and control signals define the operating scenarios. This choice is physically motivated by the role of $T_c$ as a key cold-side design and operating parameter affecting cycle efficiency, component operating states, and transient response (Zhang et al., 2022).

This paper develops a Modelica-based Helium-Xenon closed Brayton cycle model (Fig. 4) as the nonlinear transient case, while its detailed governing equations and validation are provided in our previous work (Zhang and Wang, 2024; Zhang et al., 2025a). The nonlinearity of this model mainly arises from the coupled dependence of working-fluid properties, rotating machinery characteristic curves, and loop-level thermal-hydraulic processes. The control inputs include load demand, as well as PID-based cooling-power regulation and reactor-reactivity control, leading to transient evolution governed by strongly coupled variations in both states and inputs.

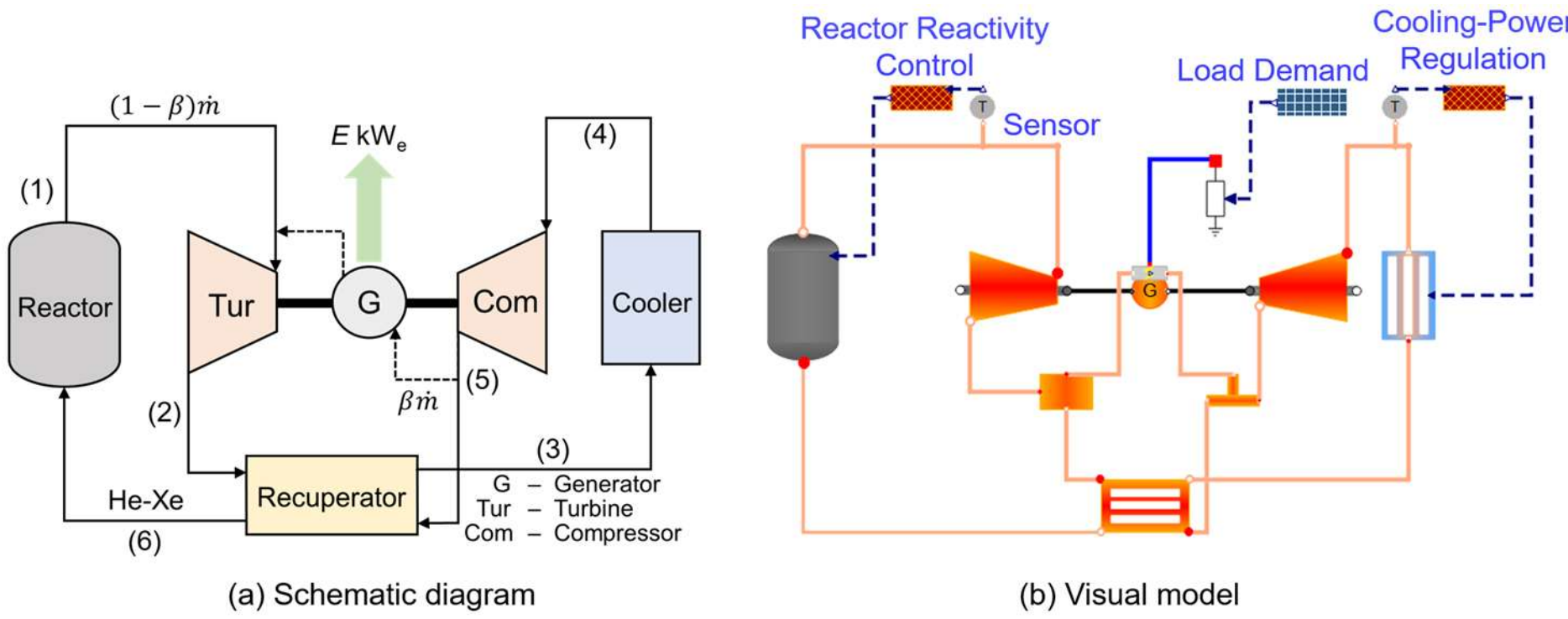

**Fig. 4.** Modelica-based visual model of the Helium-Xenon closed Brayton cycle.

Under such multi-physics coupling, a single parameter transformation cannot fully regularize all nonlinear dependencies of the reduced dynamics. To account for an important residual dependence associated with closed-loop cooling regulation, the output of the PID controller for cooling-power regulation is incorporated into $\beta$ as an auxiliary measured input $\beta^+$. Its influence is thereby represented through the input channel rather than absorbed into the state-transition operator $A$, allowing $A$ to more clearly represent the intrinsic reduced dynamics and improving prediction under the prescribed closed-loop operating scenarios.

### 2.3.3 Nonlinear Transient and Periodic System: Kármán Vortex Street

Flow past a circular cylinder is a representative nonlinear periodic system in fluid dynamics, whose wake evolves into a stable Kármán vortex street with alternating vortex shedding. This case (Fig. 5) is obtained using a Lattice Boltzmann Method (LBM, Krüger et al., 2017) solver for two-dimensional incompressible flow and checked against the open-source reference model (Luo et al., 2024) at a representative vortex-shedding condition of $Re = 640$, yielding a Strouhal number of 0.2188 close to the reference value of 0.2219, with a relative deviation of approximately 1.4%.

In this case, the physical viscosity is treated as the system parameter defining different Reynolds-number conditions, while the inlet-velocity step serves as the transient operating protocol. The characteristic geometry and the inlet-velocity protocol are kept unchanged across parameter conditions, so that the wake exhibits both transient adjustment and periodic vortex shedding. This case is therefore well suited for evaluating the proposed pROM framework for nonlinear periodic dynamics.

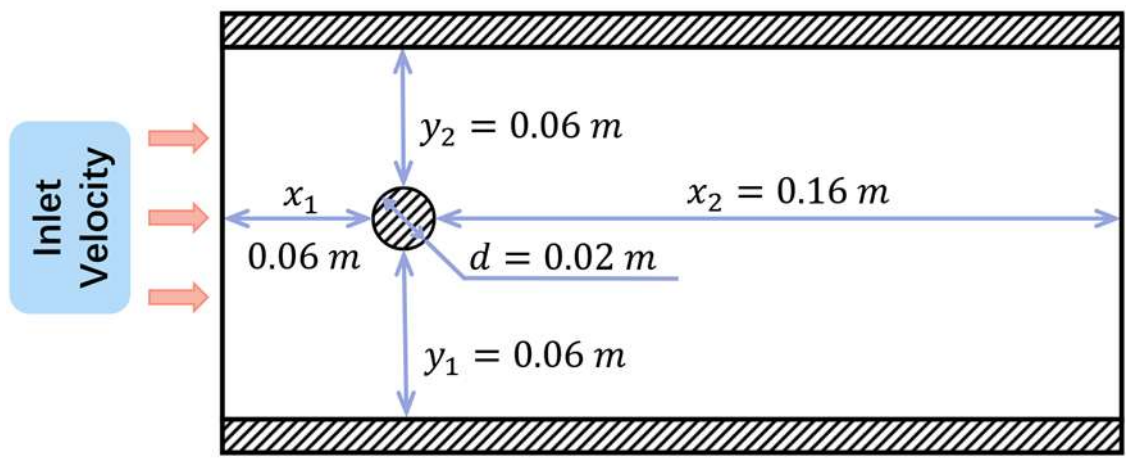

**Fig. 5.** Geometric structure of the Kármán vortex street (Luo et al., 2024).

Owing to its transient periodic nature, a nondimensional time mapping (Eq. (4)) is introduced to reduce the inconsistency of the identified spectral operator caused by frequency variation. For the Kármán vortex street, the convective time scale $\tilde{t}$, defined from the inlet velocity $U$ and cylinder diameter $D$, normalizes the basic kinematics, and the transient wake evolution under different conditions is further mapped onto a unified dynamical scale $\tau$ using the Strouhal number $St$ (Gabbai and Benaroya, 2005):

$$\tau = St\,\tilde{t} = St\left(\frac{U}{D}t\right) \tag{8}$$

Recasting the dynamics into this mapped coordinate $\tau$ yields a more consistent representation across the transient regime. Meanwhile, the amplitude envelope and $St$ are identified and extrapolated separately from the reduced linear evolution, allowing the reduced operator to focus on phase-coherent dynamics while the envelope variation is treated explicitly.

For parameter generalization, the physics-guided mapping introduced in Eq. (6) recasts the varying physical viscosity in terms of the Reynolds number ($Re$). In fluid mechanics, $Re$ is the natural dimensionless coordinate for organizing dynamic similarity and wake-regime evolution. Rather than treating viscosity as an isolated numerical input, this mapping provides a more physically meaningful axis for data-driven regression. Accordingly, the spectral properties of the wake dynamics, such as vortex-shedding frequency and amplitude evolution, vary more coherently over the $Re$ domain. This supports the Secondary DMD in extrapolating the reduced operators along a more interpretable and physically consistent parametric trajectory, thereby preserving the dominant vortex structures under unseen conditions.

## 3. Results and Discussion

### 3.1 Mechanical Transmission: Torque-Input Transients under Motor-Load Inertia Matching

For the Modelica-based mechanical transmission system, original and new scenarios under step-torque inputs are constructed together with the corresponding sampled and target parameter points (Fig. 6), parameterized by the normalized inertia ratio $J = J_m k^2/(J_i + J_l)$ to characterize motor-load inertia matching. In the original scenario, sampled points are used for model identification and additional target points are reserved for prediction. In the new scenario, which represents an unseen torque-input condition, the target responses are predicted directly using the operators identified from the sampled points of the original scenario. This setup enables a comparative evaluation of different spectral operator regression methods and the physics-guided coordinate transformation based on the torsional snapshots defined in Section 2.3.1, with time-domain predictions and error/rank comparisons reported in Figs. 7-8. The relative norm error is defined as $e_{\mathrm{rel}} = \| \hat{\mathrm{X}} - \mathrm{X}_{\mathrm{ref}} \|_F / \| \mathrm{X}_{\mathrm{ref}} \|_F$, where $\mathrm{X}_{\mathrm{ref}}$ and $\hat{\mathrm{X}}$ denote the reference and predicted time-series matrices over all evaluated variables and time instants, respectively.

The results show that the DMDc-centered pROM framework can preserve the predictive capability of reduced operators identified from the original scenario when transferred to the new scenario under different torque inputs. In interpolation (normalized $J = 0.1375$), smooth local parametric evolution allows the target operators to be accurately inferred from neighboring sampled points, so all regression methods yield similarly low errors. By contrast, in extrapolation (normalized $J = 0.5$ and 1), direct regression in the original physical $J$ coordinate leads to significant error growth, whereas the transformed coordinate $\mu' = 1/J$ reduces the dominant nonlinear parameter dependence and improves both predictive accuracy and stability, keeping relative norm errors within 1% for both interpolative and extrapolative predictions. In the mapped space, Secondary DMD outperforms conventional linear and RBF regression by modeling the structured evolution of reduced operators across the parameter domain rather than relying only on static pointwise fitting. As the dominant nonlinear dependence is weakened, the remaining parametric variation becomes lower-dimensional, making the mapped Secondary DMD less sensitive to the parametric reduced rank $r_p$. At the same time, an excessively high temporal reduced rank $r_t$ lowers predictive accuracy, suggesting that higher-order temporal modes are more sensitive to local fitting noise and therefore less robust in long-range extrapolation.

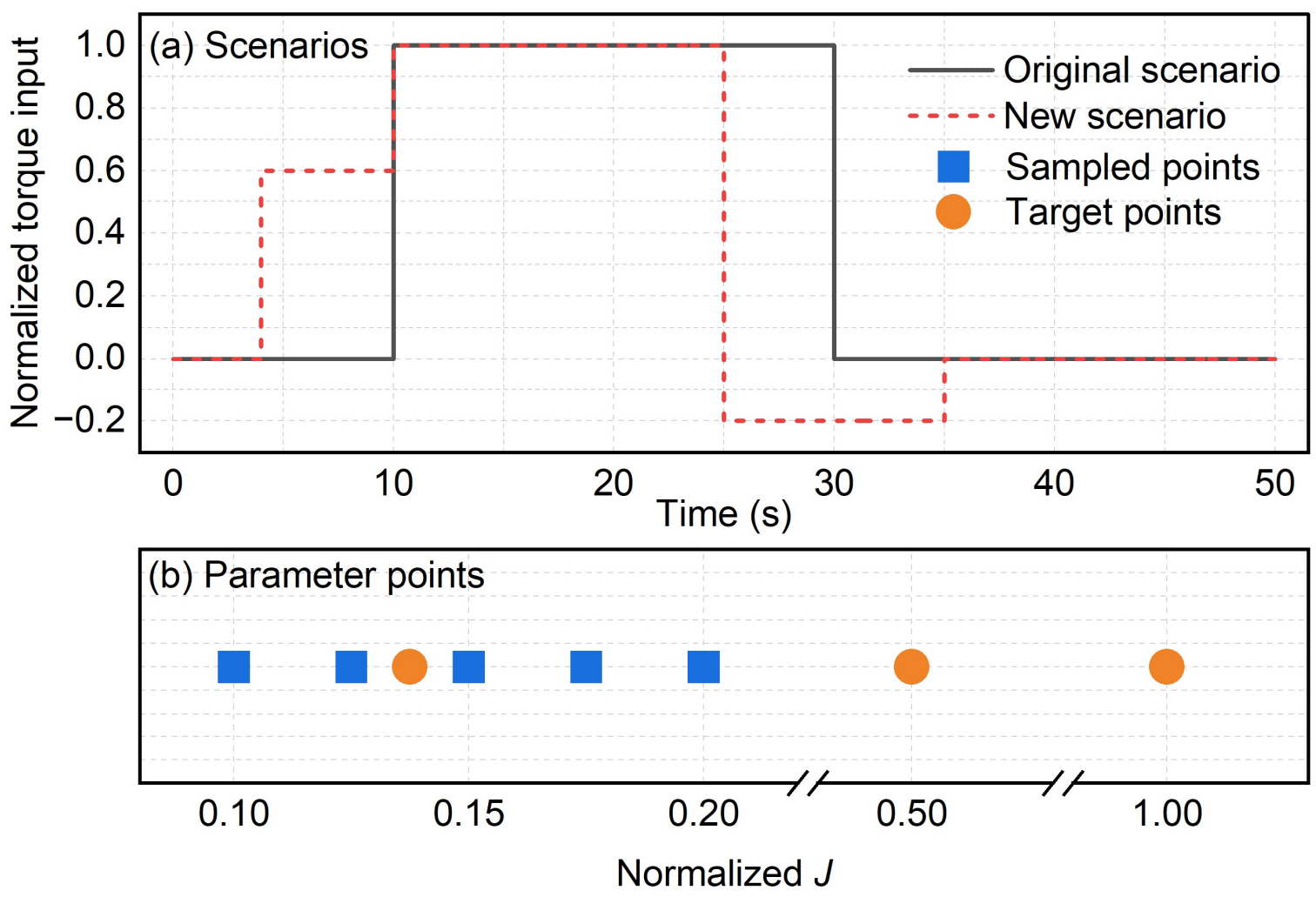

**Fig. 6.** Setup of the scenarios and parameter points.

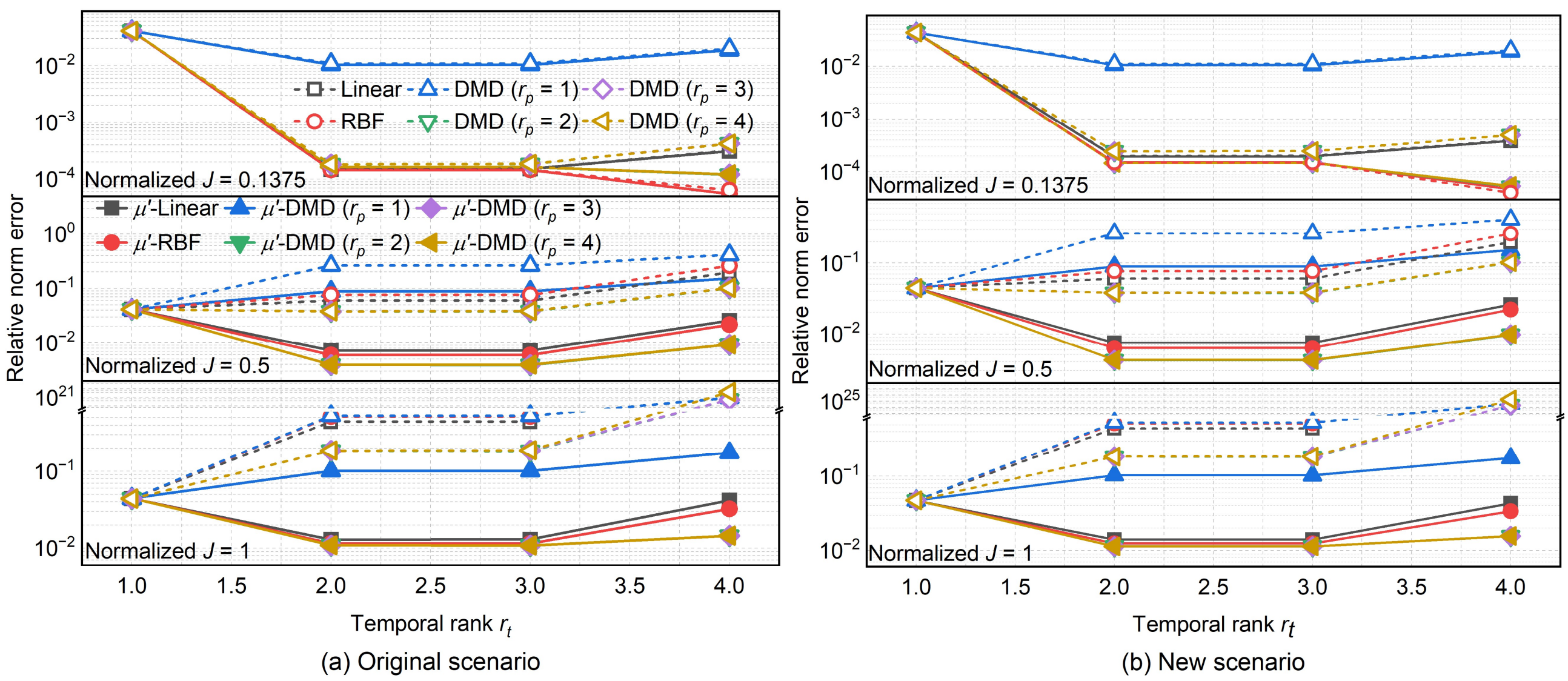

**Fig. 7.** Prediction error comparison of pROMs using different regression methods.

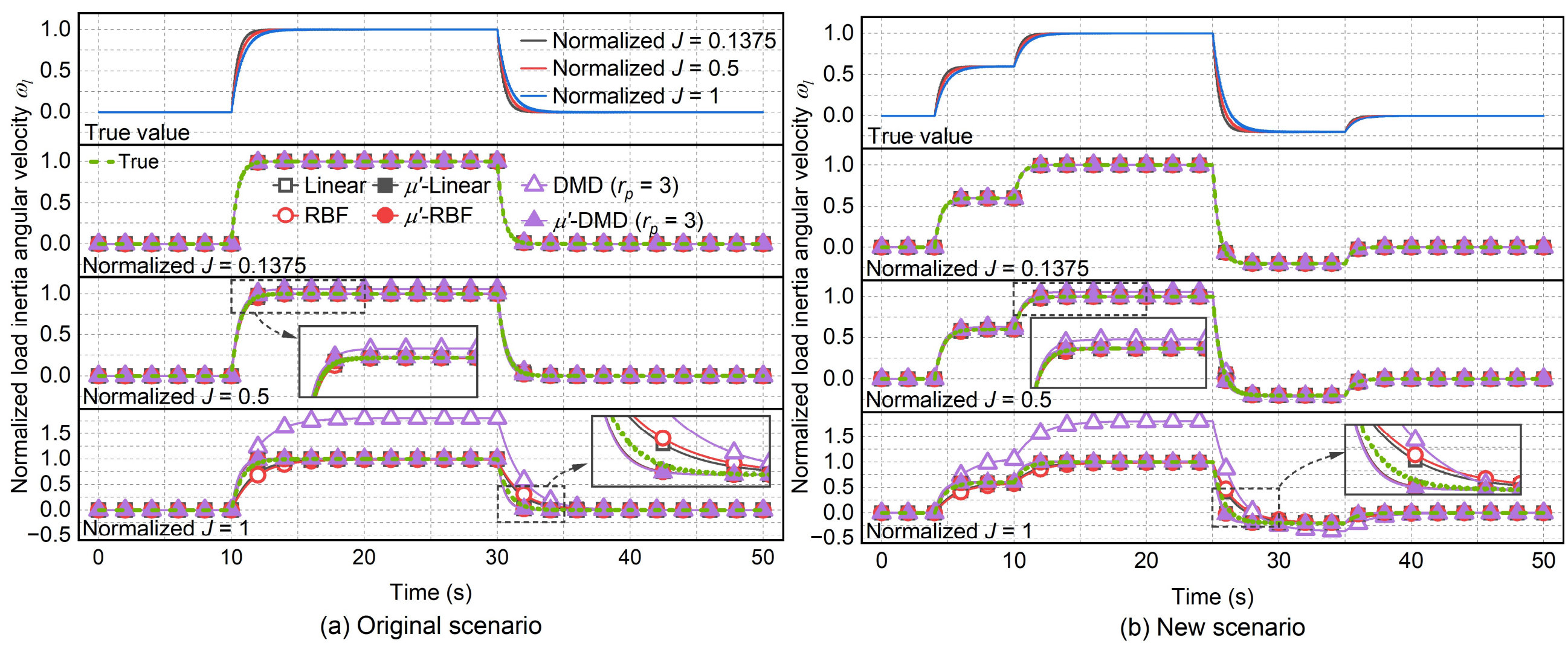

**Fig. 8.** Time-domain prediction of a representative system variable ($r_t$ = 3).

### 3.2 Closed Brayton Cycle: Load Transients under Minimum-Cycle-Temperature Variation

For the Modelica-based Helium-Xenon closed Brayton cycle, two operating scenarios are constructed, and the parameter conditions are defined by the minimum cycle temperature $T_c$ (Fig. 9). The sampled points in the original scenario are used for model identification, while target points in both the original and new scenarios are used to evaluate prediction under unseen parameter values and operating conditions. To address the strong thermohydraulic coupling, the augmented input vector further includes the measured auxiliary inputs $\beta^+$ defined in Section 2.3.2.

Under the augmented formulation, the relative norm errors in Fig. 10 highlight the critical role of baseline regularization ($\Delta$, Eq. (3)). In interpolation ($T_c = 374.5$ K), all regression methods yield similarly low errors because local parametric evolution remains smooth. By contrast, in extrapolation ($T_c = 390$ and $400$ K), the $\Delta$-regularized formulation substantially outperforms the uncentered formulation. By removing steady-state offsets, the regularized operators can focus on intrinsic transient dynamics, maintaining relative norm errors consistently below 1%. Fig. 10 also shows that higher temporal ranks ($r_t \geq 3$) are highly sensitive to local fitting noise, making error growth difficult to control and indicating the need for strict low-rank truncation, which is consistent with the observations in Section 3.1.

Fig. 11 further shows the time-resolved relative errors of all predicted system variables at the farthest extrapolation point $T_c = 400$ K, with pointwise errors remaining within ±3% in both the original and new scenarios. Fig. 12 presents the corresponding time-domain comparison for the key variable shaft speed, showing good agreement between the predicted and reference responses.

With the $\Delta$-regularized formulation fixed, Fig. 13 further illustrates the effect of $\beta^+$. The comparison between the full framework and the formulation without $\beta^+$ shows that the latter still captures the main transient trend, but the inclusion of $\beta^+$ leads to much closer agreement with the reference response. This result highlights the importance of measured auxiliary inputs for improving predictive accuracy in strongly coupled thermohydraulic systems.

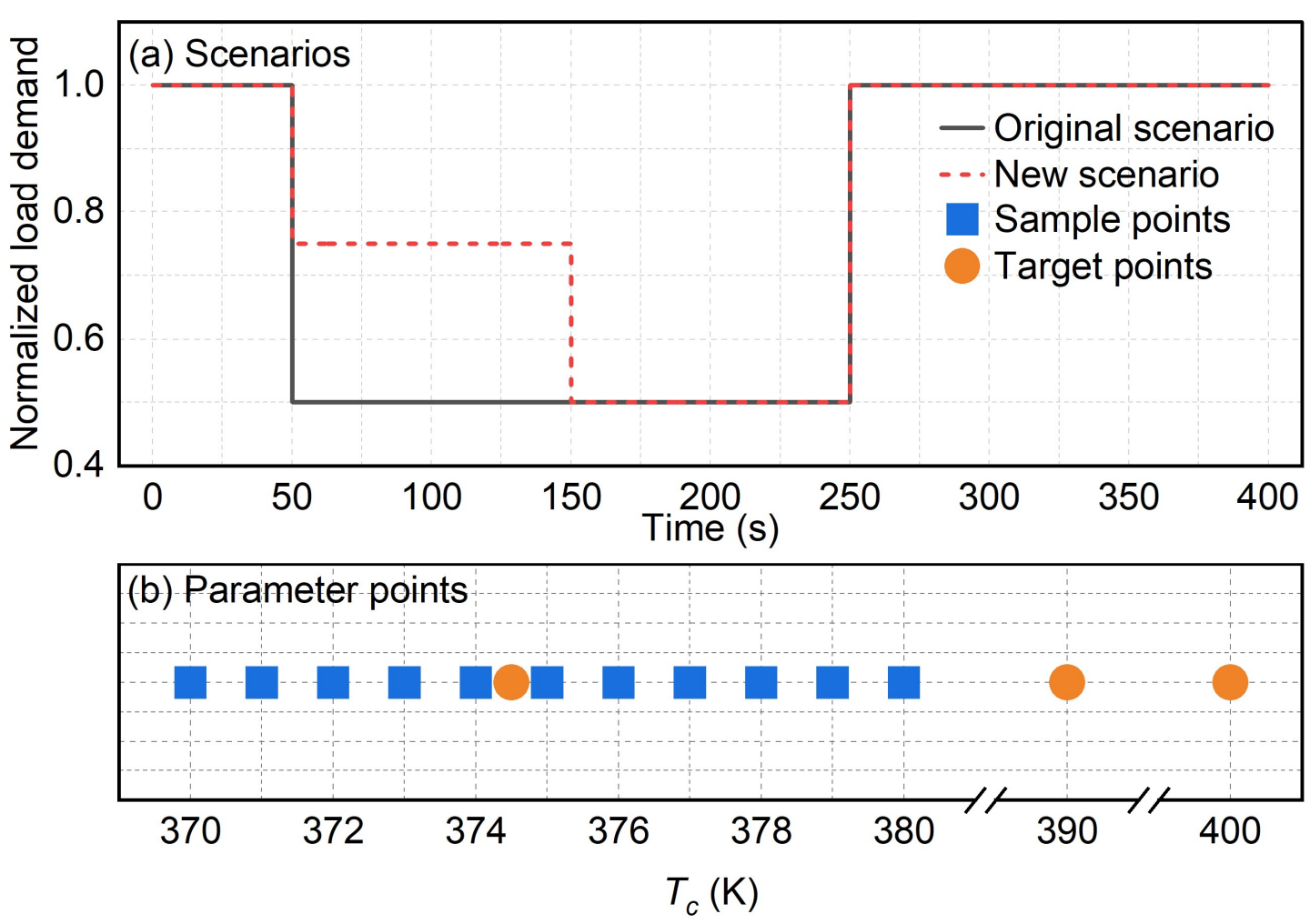

**Fig. 9.** Setup of the scenarios and parameter points.

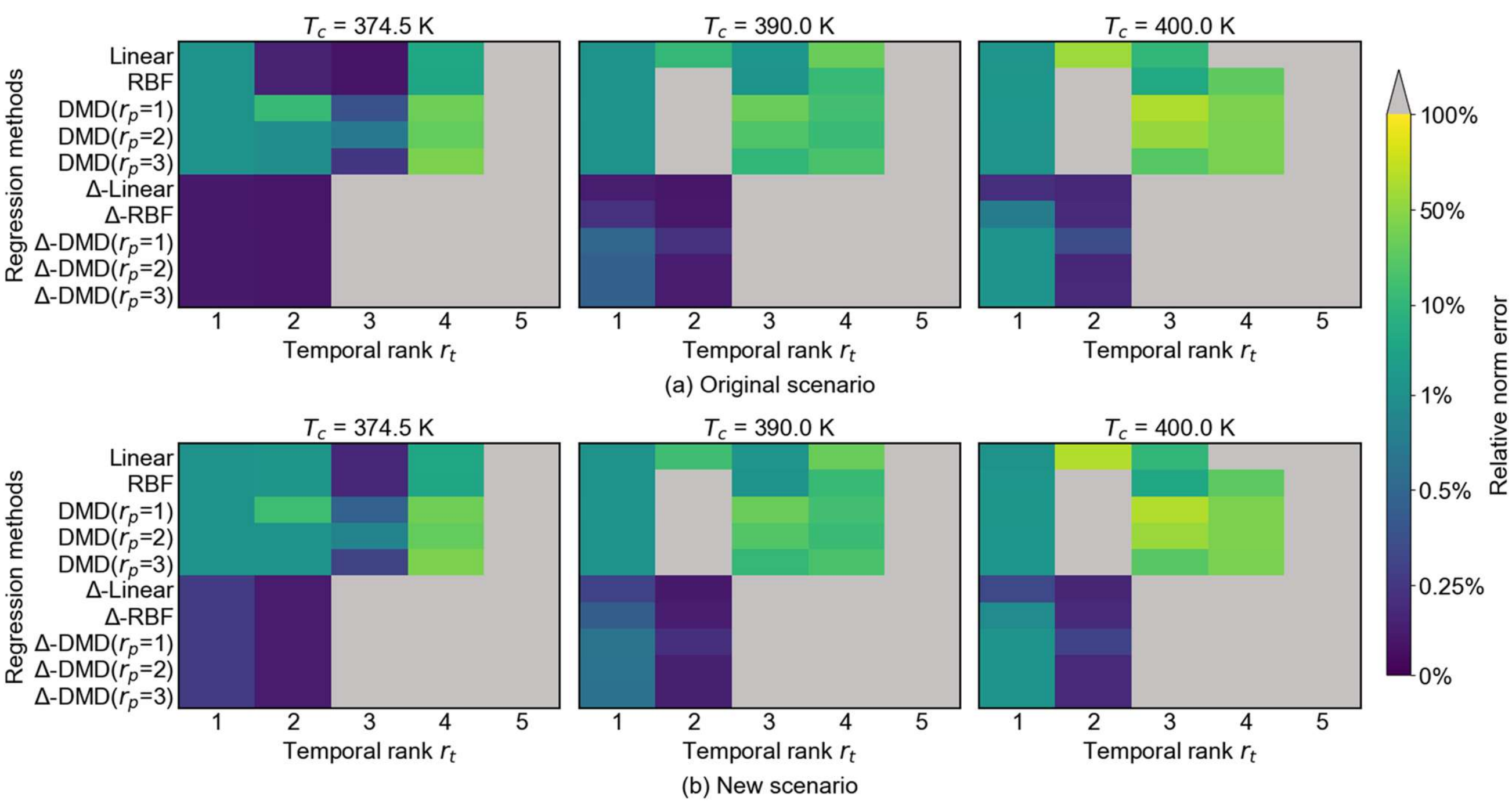

**Fig. 10.** Prediction error comparison of pROMs using different regression methods.

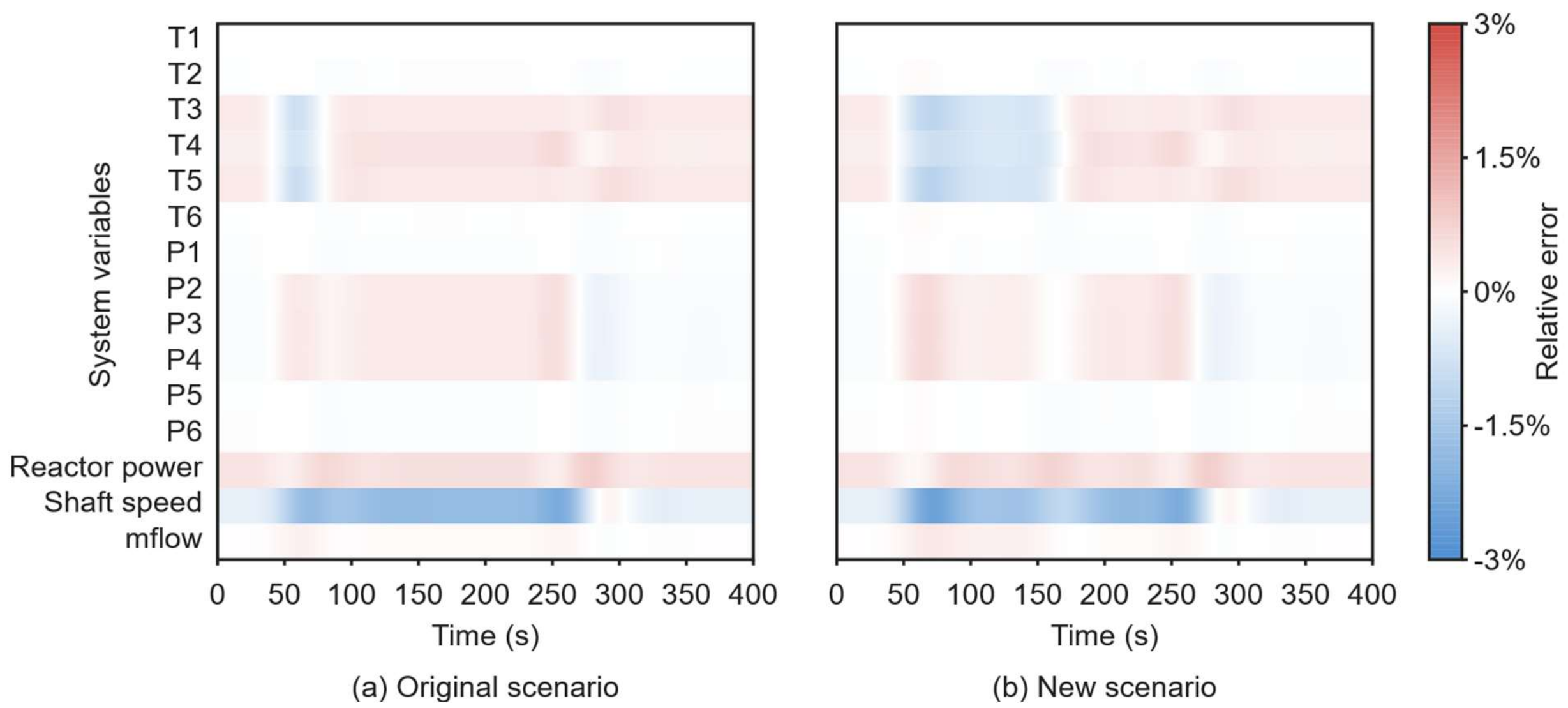

**Fig. 11.** Time-resolved relative errors of all predicted system variables at $T_c = 400$ K ($r_{t,p} = 2$).

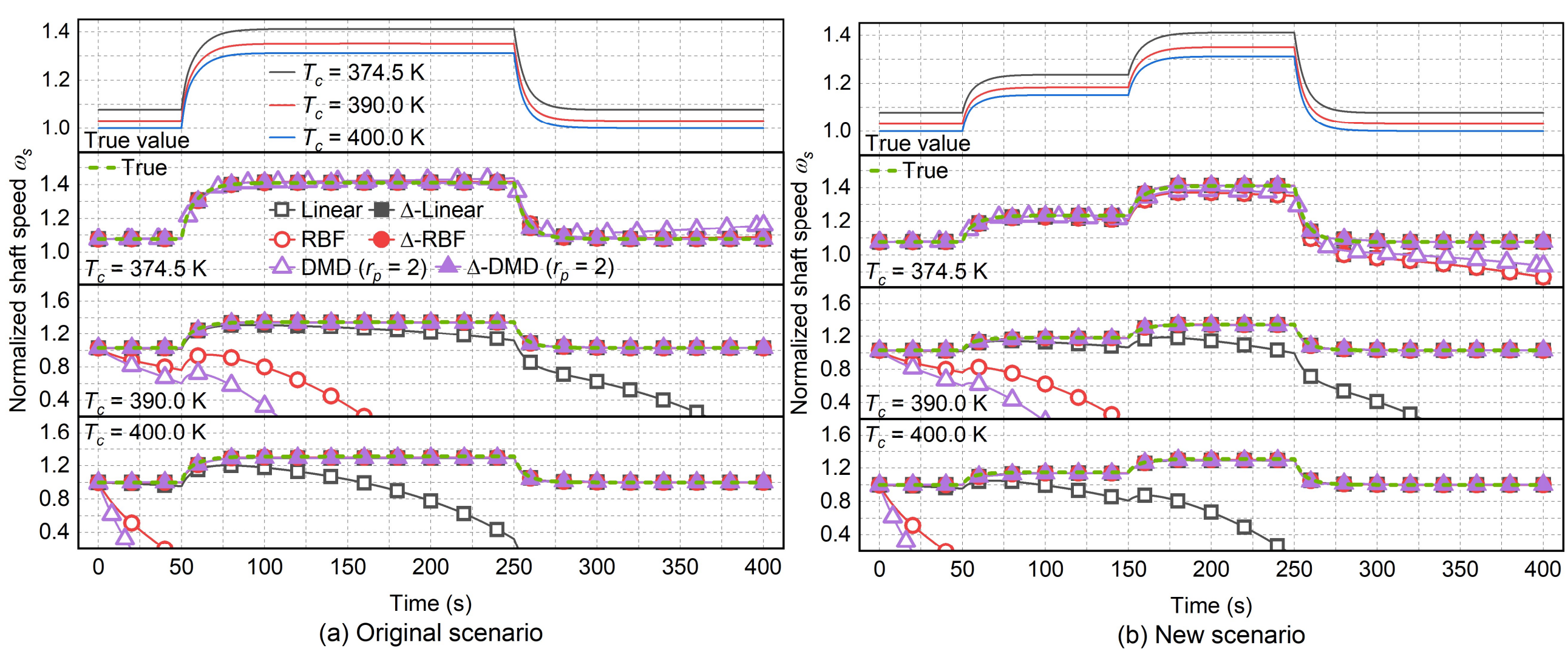

**Fig. 12.** Time-domain prediction of a representative system variable ($r_t = 2$).

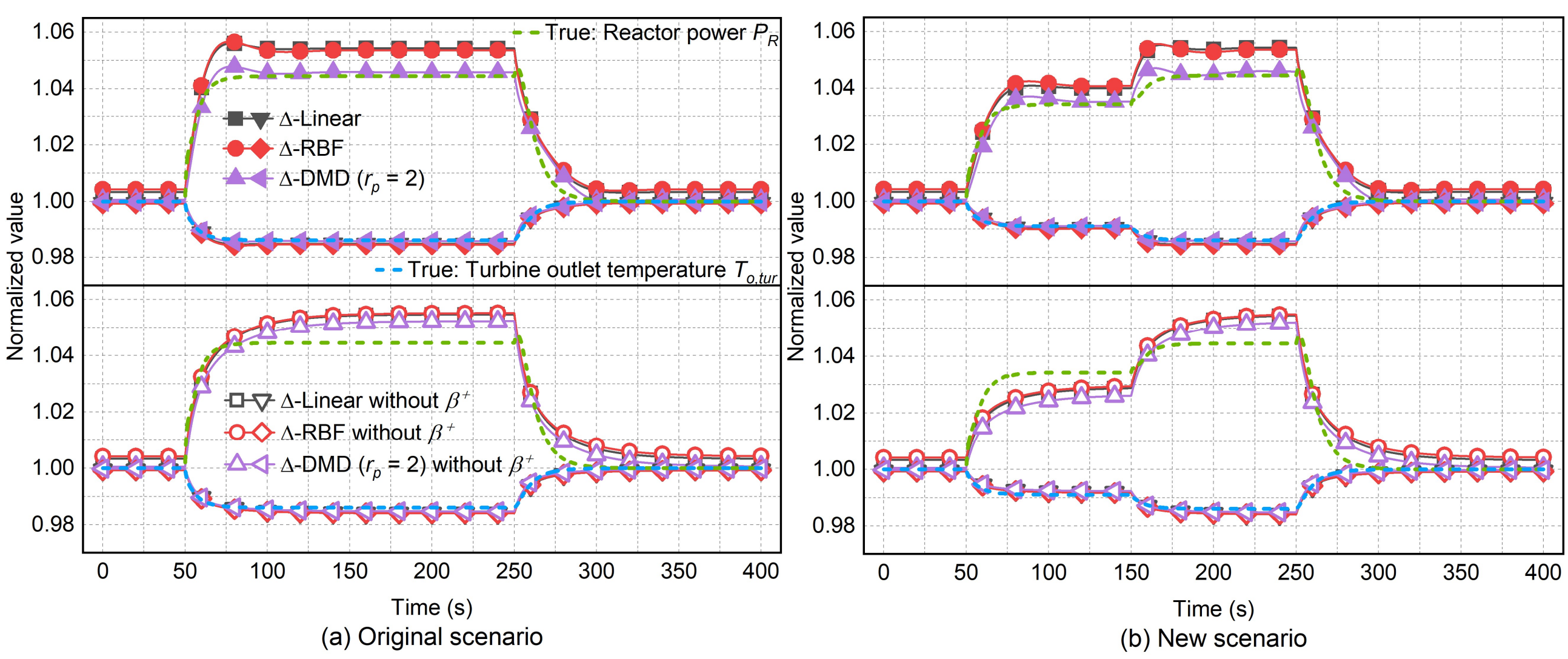

**Fig. 13.** Effect of measured auxiliary inputs $\beta^+$ ($r_t = 2$).

### 3.3 Kármán Vortex Street: Inlet-Velocity Transients under Fluid-Viscosity Variation

In the Kármán vortex-street case, the reduced-order model is constructed from full-field velocity snapshots of the cylinder wake under step changes in inlet velocity, while three velocity probes are introduced only as supplementary local diagnostics of the reconstructed transient dynamics (Fig. 14). The physical viscosity is transformed into the Reynolds number $Re$ as the parameter coordinate for regression. For vortex streets under constant inlet velocity, standard DMD can provide sufficient temporal identification and reliable prediction (Riva et al., 2026). Under transient inlet-velocity conditions, however, direct DMD is no longer sufficient, and standard DMDc shows degraded performance because the vortex-shedding time scale varies with the changing flow condition. After introducing the nondimensional time mapping onto $\tau$ (Eq. (8)), this temporal inconsistency is reduced, allowing $\tau$-DMDc to provide more accurate temporal identification (Fig. 15). The parameter-direction Secondary DMD regression further gives lower relative norm errors than linear and RBF regressions (Fig. 16). The predicted spatiotemporal vortex-shedding fields remain consistent with the reference evolution (Fig. 17), and the velocity magnitudes at selected probes are also accurately reproduced (Fig. 18). These results confirm that the proposed framework captures the transient variations of the vortex street, while the extrapolative predictions preserve the dominant shedding structures despite some loss of local accuracy.

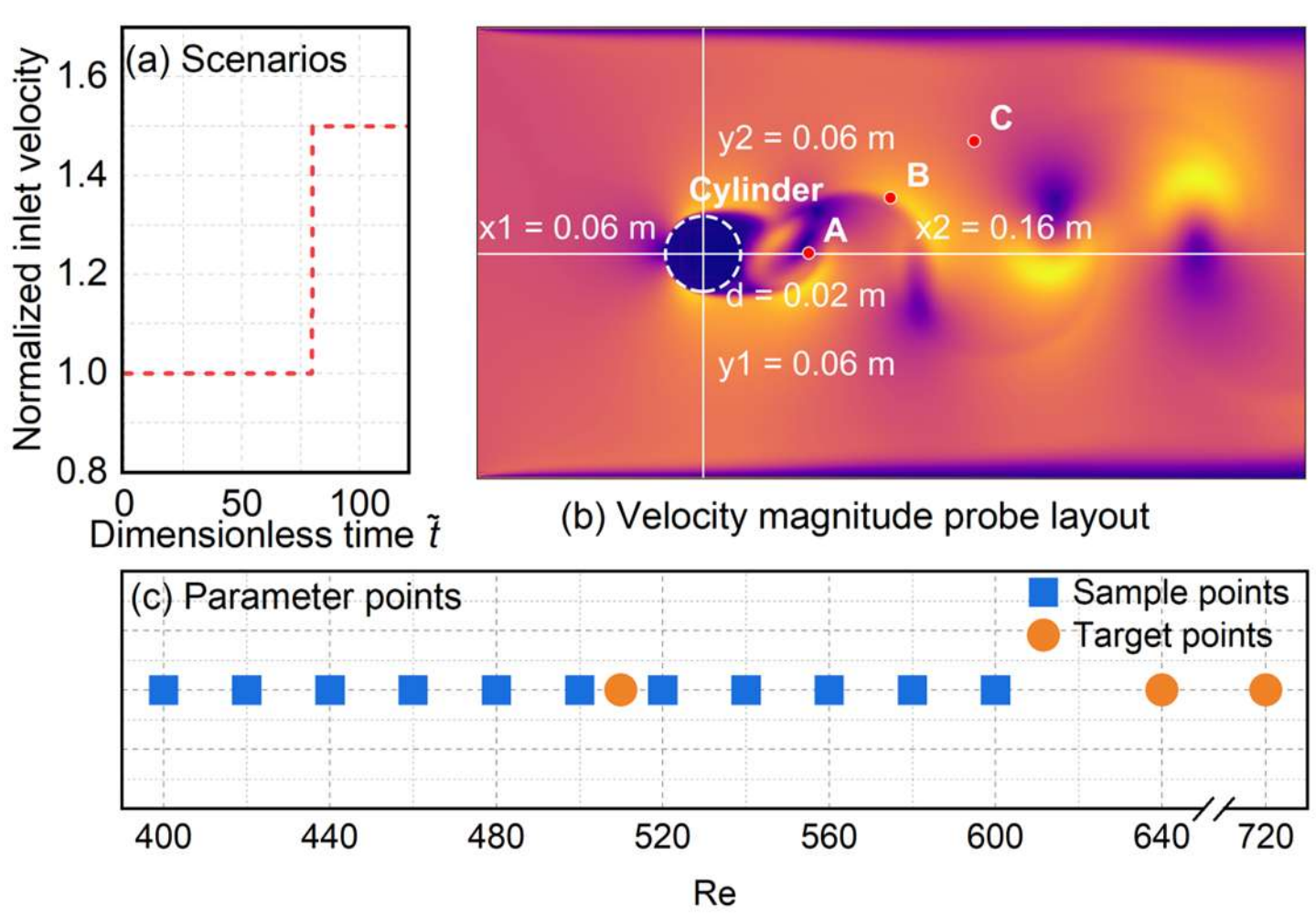

**Fig. 14.** Setup of the scenarios and parameter points.

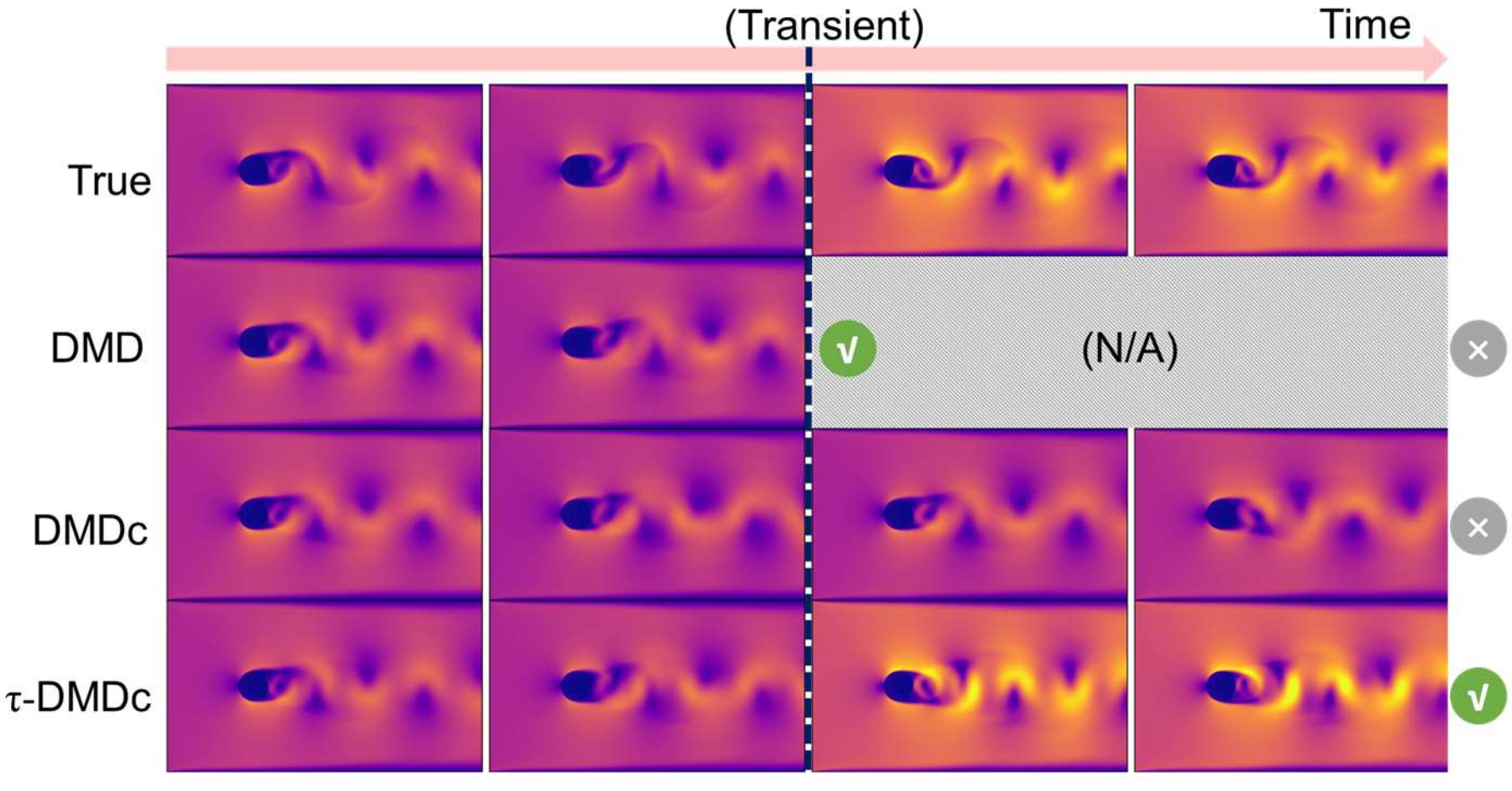

**Fig. 15.** Effect of nondimensional time mapping on transient periodic prediction.

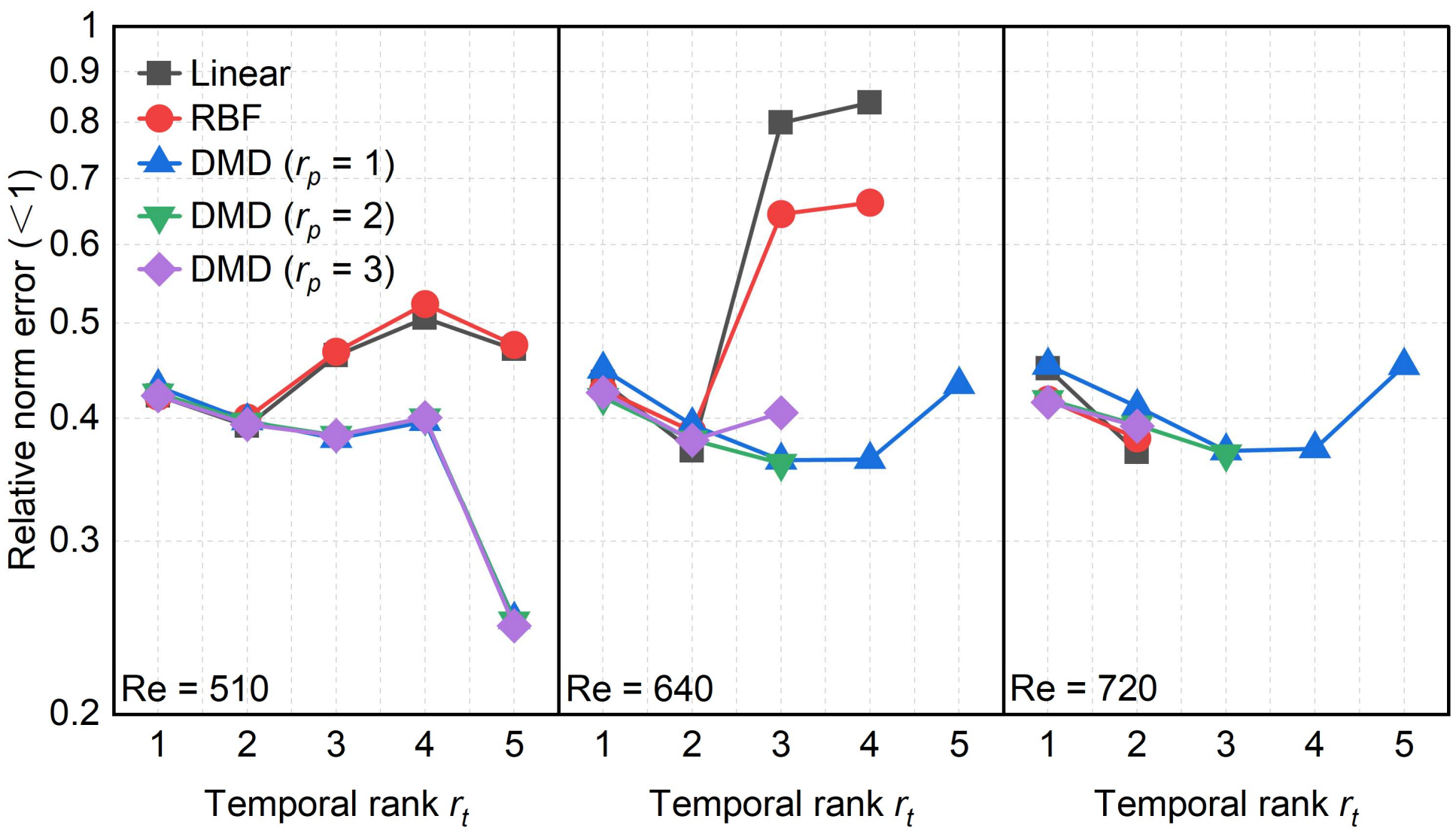

**Fig. 16.** Prediction error comparison of pROMs using different regression methods.

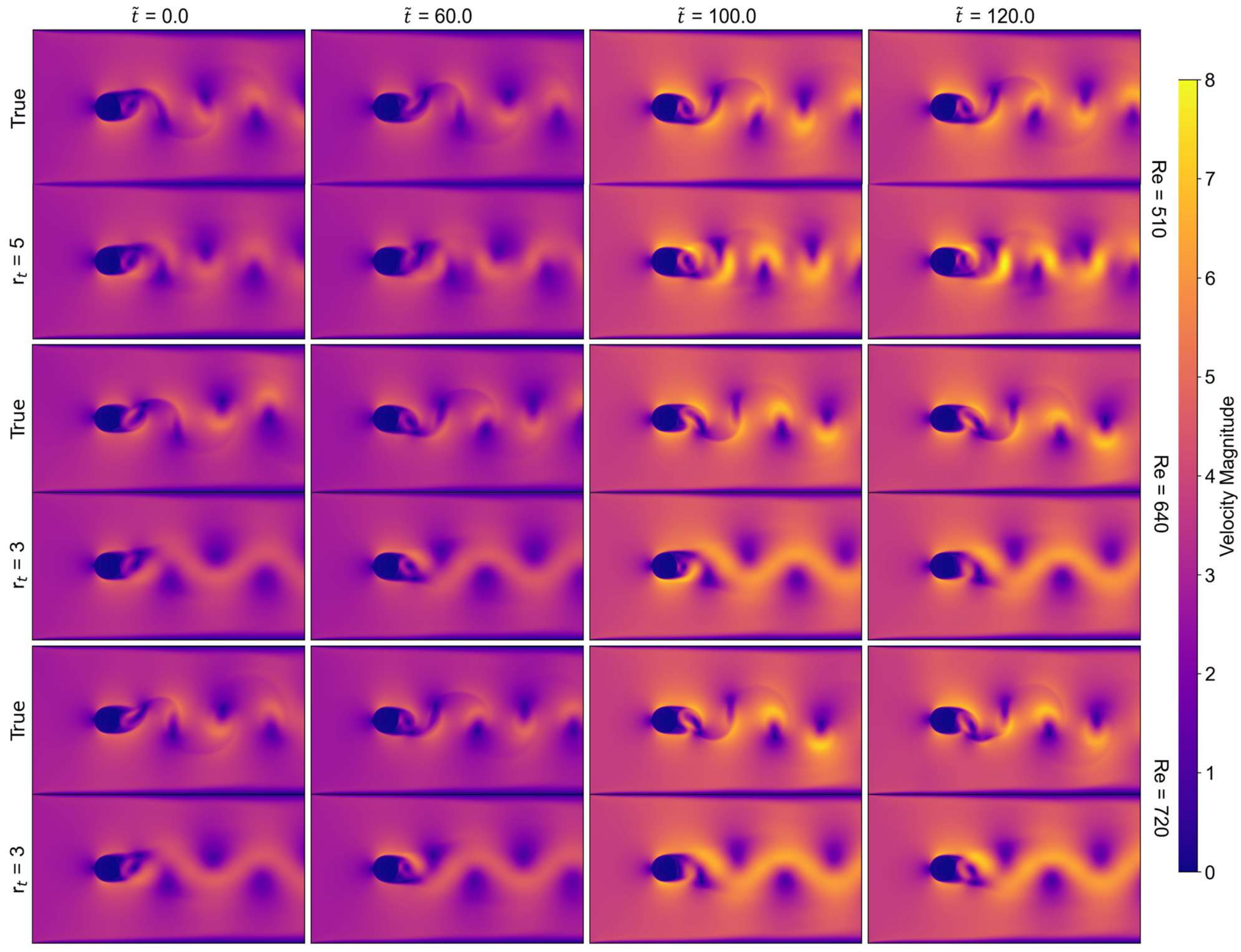

**Fig. 17.** Time-domain prediction of the Kármán vortex street (Secondary DMD, $r_p = 2$).

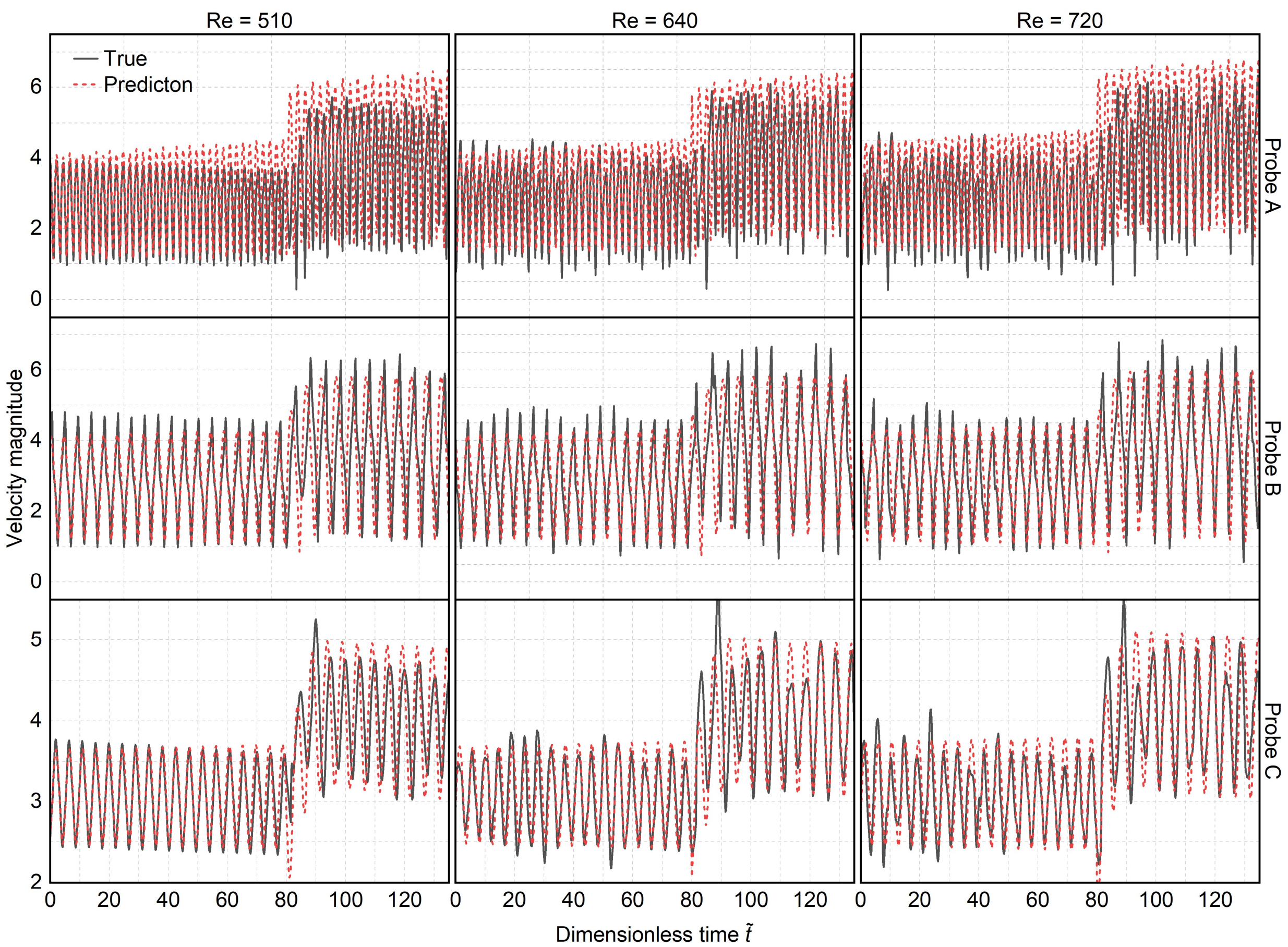

**Fig. 18.** Predictions of velocity magnitude at selected probes.

### 3.4 Further Discussion

Building on the predictive capability demonstrated across linear, nonlinear, and periodic archetypes, this section further discusses the data requirements of the proposed physics-guided spectral pROM framework, its comparison with machine-learning surrogates, and confidence assessment for extrapolative prediction.

#### 3.4.1 Sampling and Noise Considerations

The present work focuses on parametric transfer of reduced spectral operators rather than on exhaustive sample-size optimization. The sampled conditions are therefore selected to provide a stable local spectral basis for each representative case, while the influence of sample size and measurement noise on DMD/DMDc identification has been systematically discussed in our previous study (Zhang et al., 2025a). In practical applications, insufficient samples or noisy snapshots may degrade spectral alignment and should be addressed through additional sampling or denoising procedures such as Kalman-filter-assisted preprocessing.

### 3.4.2 Comparison with LSTM-Based Machine-Learning Prediction

Compared with the proposed reduced-order operator-based framework, black-box machine-learning surrogates typically offer weaker interpretability and rely more heavily on the coverage of the training dataset, which may limit their reliability for prediction outside the sampled parameter and operating-condition space. To provide a direct comparison, an LSTM-based time-series prediction model is trained and tested on the Modelica-based closed Brayton cycle under the same interpolation/extrapolation targets and operating scenarios. The reported LSTM results are obtained after hyperparameter tuning, and the best-performing configuration is used for comparison. As shown in Fig. 19, the proposed pROM framework achieves more accurate predictions than the LSTM-based model, especially for extrapolation and new operating conditions. This comparison further supports the advantage of physics-guided reduced spectral-operator transfer over purely data-driven sequence learning when only limited sampled parameter conditions are available.

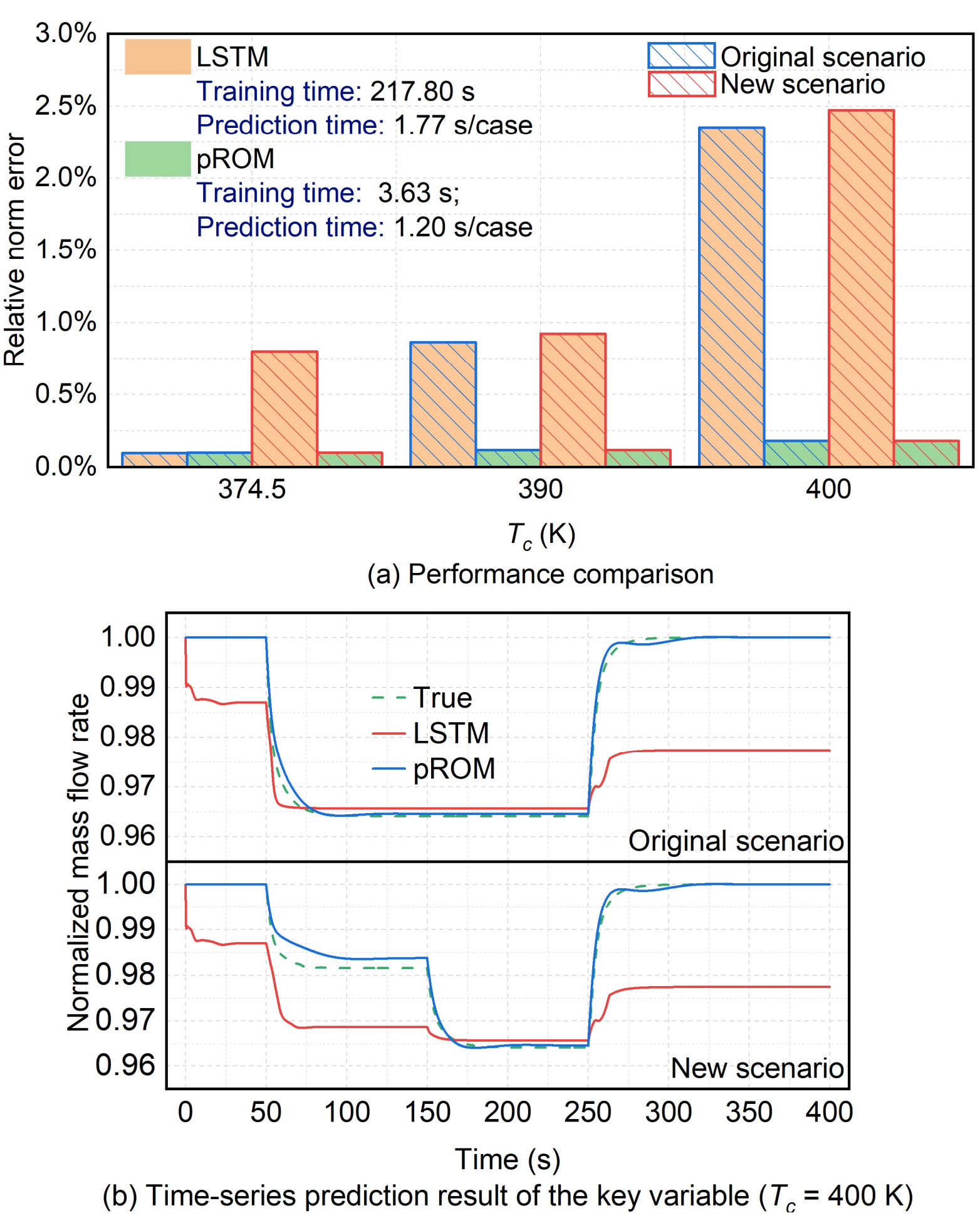

**Fig. 19.** Comparison between pROM and LSTM.

### 3.4.3 Confidence Assessment for Extrapolation

Prediction reliability is central to parametric extrapolation. As shown in the results, extrapolation errors generally increase as the target parameter moves farther from the sampled domain, while the credibility of such predictions depends not only on regression accuracy, but also on whether the target condition remains within a physical regime where the dominant modal structure is preserved. For example, in the mechanical transmission case, the normalized $J$, characterizing the system's inertia matching, ranges from 0 to 1 and thus provides a physically meaningful reference interval for modal consistency. When such boundaries cannot be explicitly identified, a practical alternative is to reserve independent test points, map relative norm errors over the parameter axis, and construct an auxiliary error-prediction surrogate to indicate extrapolation confidence, following a reliability-assessment rationale also used in applicability-domain analysis based on independent validation errors (Sutton et al., 2020). Taking the Helium-Xenon closed Brayton cycle case (Section 3.2) as an example, a few additional extrapolation points can be selected as test points, and fitting their errors provides an effective estimate of the expected error level for farther extrapolation target points (Fig. 20), thereby supporting dataset updates when the estimated error exceeds a predefined tolerance.

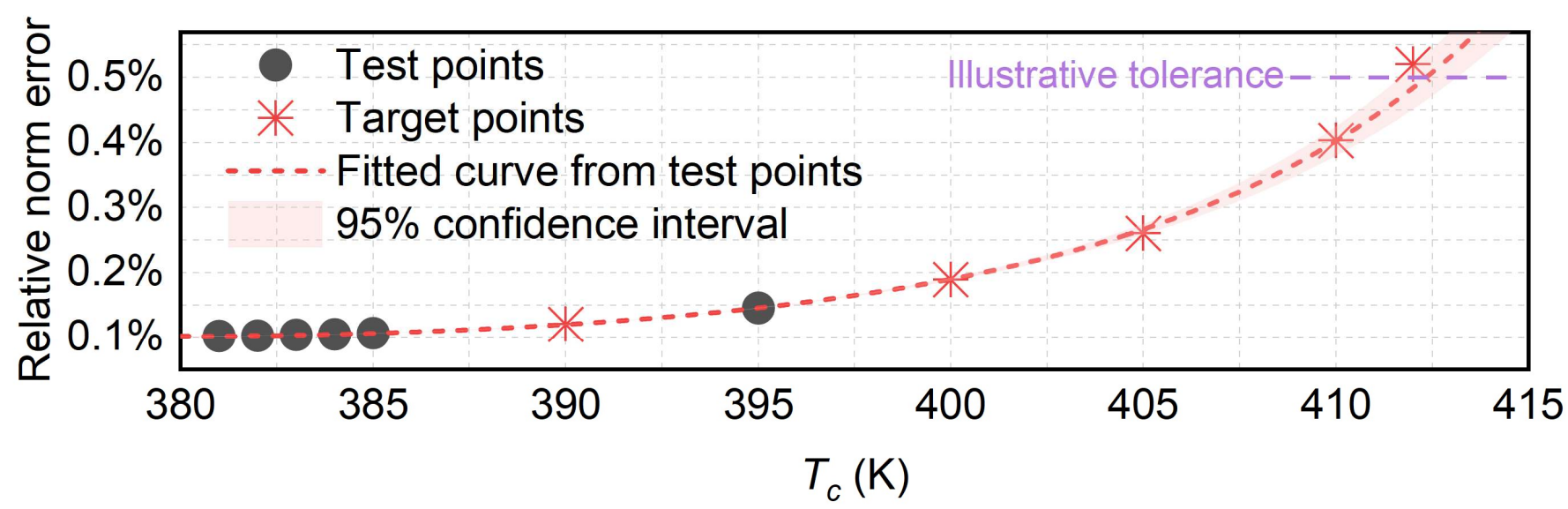

**Fig. 20.** Extrapolation confidence assessment for the Helium-Xenon closed Brayton cycle (Secondary DMD, $r_{t,p} = 2$).

## 4. Conclusions

This paper establishes a DMDc-centered physics-guided spectral pROM framework for interpretable system identification and extrapolative transient prediction of controlled dynamical systems. By separating intrinsic reduced dynamics from external control effects and transferring coherent reduced spectral quantities across parameter conditions, the framework enables prediction at unseen parameter values and under new operating conditions using limited sampled data. The main conclusions are as follows:

(1) The proposed framework identifies spectrally coherent reduced-order dynamics across linear, nonlinear transient, and nonlinear periodic systems. Physics-guided coordinate transformations, baseline regularization, and nondimensional time mapping further improve predictive robustness by reducing nonlinear parameter dependence, steady-state offsets, and temporal-scale inconsistency.

(2) The framework achieves strong extrapolative prediction performance across three representative dynamical systems. For the mechanical transmission and Helium-Xenon closed Brayton cycle cases, relative norm errors remain below 1% in both interpolation and extrapolation under unseen parameter values and new operating conditions. In particular, the Brayton-cycle case demonstrates that the proposed framework can predict system-level multivariable transients involving thermohydraulic coupling, rotating-machinery characteristics, and closed-loop regulation under unseen cold-side temperature conditions and new load-demand scenarios. For the Kármán vortex-street case, the dominant vortex-shedding structures and main spatiotemporal evolution are preserved despite phase-induced local errors.

(3) Secondary DMD provides a clear advantage over linear and RBF regressions for parameter-domain operator prediction. This advantage is especially important in extrapolation, where target operators are not bracketed by neighboring sampled conditions and prediction must rely on learned spectral evolution beyond the observed parameter interval. The comparison with an LSTM-based machine-learning surrogate further supports the improved extrapolative robustness of the proposed pROM under limited sampled parameter conditions.

(4) The reduced-rank analysis shows that dominant low-order modes provide the most robust basis for extrapolative prediction. The optimal reduced rank remains largely consistent across extrapolation targets, reducing model-selection ambiguity and supporting nearby short-range extrapolation tests for prediction-confidence assessment. In contrast, excessive temporal ranks amplify local fitting noise and degrade predictive robustness.

Overall, the proposed framework offers an efficient and interpretable solution for extrapolative transient prediction under unseen parameter values and new operating conditions using limited sampled conditions. Future work will explore neural-network-assisted coordinate mapping and multi-parameter generalization.

**Acknowledgements**

This work was supported by the National Natural Science Foundation of China (Grant No. 12375166).